\documentclass[a4paper,10pt]{article}
\usepackage{amsmath,amscd,amssymb,french}
\newcommand{\qed}{\hfill$\square$}
\begin{document}
\begin{center}
{\Large\textbf{
STRUCTURES DE POISSON SUR LES VARI\'ET\'ES\\
ALG\'EBRIQUES DE DIMENSION TROIS}}
\end{center}
\vspace{0.5cm}
\centerline{St\'ephane DRUEL}
\begin{center}
DMI-\'Ecole Normale Sup\'erieure\\
45 rue d'Ulm\\
75005 PARIS\\
e-mail: \texttt{druel@clipper.ens.fr}
\end{center}
\vspace{1cm}
$\centerline{{\Large\textbf{\S0 Introduction}}}$\\ 
\newline
\indent Soit $X$ une vari\'et\'e alg\'ebrique lisse sur le corps
$\mathbb{C}$ des nombres complexes. Une \emph{structure de Poisson}
sur $X$ est la donn\'ee d'une structure d'alg\`ebre de Lie sur la faisceau
structural $\mathcal{O}_{X}$ de $X$, qui soit une d\'erivation en
chacune des variables. Une telle structure est d\'efinie par un tenseur
antisym\'etrique non nul $\sigma\in
H^{0}(X,\overset{2}{\wedge}\mathcal{T}_{X})$, \emph{le bivecteur de Poisson},
duquel on d\'eduit une fl\`eche $\mathcal{O}_{X}$-lin\'eaire
$\Omega_{X}^{1}\longrightarrow\mathcal{T}_{X}$. Le \emph{rang} de la
structure en $x\in X$ est, par d\'efinition, le rang de l'application
$\Omega_{X}^{1}\otimes{k(x)}\longrightarrow\mathcal{T}_{X}\otimes{k(x)}.$
La structure est dite \emph{r\'eguli\`ere} lorsque son rang en tout point
est maximal. Elle est dite \emph{quasi-r\'eguli\`ere} lorsqu'elle
est r\'eguli\`ere sauf en un nombre fini de points. Notons que le rang
d'une structure de Poisson en un point est pair 
puisqu'une telle structure est par d\'efinition antisym\'etrique. Par exemple,
une vari\'et\'e ab\'elienne $A$ de dimension au moins deux poss\`ede de nombreuses
structures de Poisson r\'eguli\`eres. En effet, toute forme multilin\'eaire
altern\'ee de rang maximal sur l'espace tangent \`a l'origine fournit
une telle structure.\\
\indent Une surface projective
admettant une structure de Poisson non triviale est une surface $K3$,
une surface ab\'elienne ou une surface r\'egl\'ee; une telle structure
est alors d\'efinie par la seule donn\'ee d'une
section du fibr\'e anticanonique.\\
\indent Ce travail est consacr\'e \`a l'\'etude des vari\'et\'es alg\'ebriques
projectives lisses de dimension 3, admettant une structure de
Poisson non triviale. Une telle structure \'etant d\'efinie par un tenseur
antisym\'etrique d'ordre 3, on s'attend \`a ce que ladite structure ne
s'annule qu'en un nombre fini de points, seul cas que nous
\'etudierons. Le r\'esultat principal est le th\'eor\`eme:\\
\newline
\indent\textbf{Th\'eor\`eme 0.1 }\textit{ Soit $X$ une vari\'et\'e projective
  lisse de dimension 3.\\
\indent Alors $X$  admet une structure de Poisson
  quasi-r\'eguli\`ere si et seulement si  $X$ appartient \`a l'une des 4
  familles suivantes:\\ 
\indent (1) $X$ est une vari\'et\'e ab\'elienne,\\
\indent (2) $X$ est un fibr\'e plat en droites projectives sur une
surface ab\'elienne,\\
\indent (3) $X=(C\times A)/G$ o\`u $C$ est une courbe,
$A$ est une surface ab\'elienne et
$G\subset\text{Aut}(C)$ un groupe fini op\'erant librement sur $C\times
A$ par la formule: 
$$g.(c,a)=(g.c,t_{g}(c,a)+u_{g}(a)),\,g\in G,\,c\in C,\,a\in A,$$
\noindent o\`u $u_{g}$ 
un automorphisme de groupes de $A$ respectant la structure
 symplectique et $t_{g}$
une fonction r\'eguli\`ere sur $C\times A$ \`a valeurs dans $A$,\\
\indent (4) $X=(C\times S)/G$ o\`u $C$ est une courbe,
 $S$ est une surface $K3$ et $G$ un groupe fini op\'erant librement sur
 $C$, op\'erant sur $S$ en respectant la structure symplectique et sur
 le produit $C\times S$ par le produit de ses actions sur chacun des
 facteurs.}\\
\newline
\indent La d\'emonstration de ce r\'esultat repose sur un th\'eor\`eme de
 structure de S.Mori ([Mo]) d\'ecrivant les contractions extr\'emales d'une
 vari\'et\'e projective lisse non minimale de dimension trois. Ainsi, nous
 montrons qu'une vari\'et\'e projective non minimale de dimension trois
 admettant une structure de Poisson quasi-r\'eguli\`ere est un fibr\'e en
 coniques, ce qui fait l'objet du paragraphe 2, le paragraphe 1
 donnant les propri\'et\'es essentielles de ces structures en dimension
 trois. Dans les paragraphes 3 et 4, nous \'etudions le cas o\`u la vari\'et\'e
 est minimale en \'etudiant respectivement le morphisme d'Albanese et le
 morphisme canonique.\\
\newline
\indent\textbf{Remerciements} Je tiens \`a exprimer ici tout mes
 remerciements \`a Arnaud Beauville pour son aide au cours de la
 pr\'eparation de ce travail.\\
\vspace{1cm}\\
\centerline{{\Large\textbf{\S1 Propri\'et\'es}}}
$\ $
\newline
\indent Une vari\'et\'e (alg\'ebrique) d\'esignera un sch\'ema int\`egre,
s\'epar\'e et de type fini sur le corps
$\mathbb{C}$ des nombres complexes. Nous identifierons, sans le
 mentionner, un sch\'ema et l'espace analytique 
complexe qui lui est associ\'e et nous supposerons toujours qu'une
vari\'et\'e est lisse, sauf mention du contraire.\\
\newline
\indent Soit $X$ une vari\'et\'e alg\'ebrique de dimension 3 et supposons
 que $X$ admette une structure de Poisson r\'eguli\`ere. Alors, le
 bivecteur de Poisson $\sigma\in
 H^{0}(X,\overset{2}{\wedge}\mathcal{T}_{X})$ fournit une application
 $\mathcal{O}_{X}\hookrightarrow\overset{2}{\wedge}\mathcal{T}_{X}$,
 et, compte tenu de l'isomorphisme canonique
 ${\wedge}\mathcal{T}_{X}\otimes\omega_{X}\cong\Omega_{X}^{1}$, on en d\'eduit une injection de fibr\'es vectoriels $\omega_{X}\hookrightarrow\Omega_{X}^{1}$. R\'eciproquement, une telle injection de fibr\'es vectoriels d\'efinit une structure de Poisson r\'eguli\`ere d\`es que l'identit\'e de Jacobi est satisfaite.
 Un calcul en coordonn\'ees locales permet alors de prouver le\\
\newline 
\indent\textbf{Lemme 1.1 }\textit{Soit $X$ une vari\'et\'e
  alg\'ebrique de dimension 3
admettant une structure de Poisson r\'eguli\`ere (de rang 2).\\
\indent (1) On a une suite exacte courte:\\
$$0\longrightarrow\omega_{X}\longrightarrow\Omega_{X}^{1}\longrightarrow\mathcal{T}_{X}\longrightarrow\omega_{X}^{-1}\longrightarrow0.$$}\\
\noindent\indent\textit{(2) Notons}
$\mathcal{F}=\text{Ker}(\mathcal{T}_{X}\longrightarrow\omega_{X}^{-1})$.
\textit{Alors $\mathcal{F}$ est un fibr\'e vectoriel de rang 2,
  int\'egrable, c'est \`a dire stable par le crochet de Lie naturel sur
$\mathcal{T}_{X}$.\\
\indent (3) Une section partout non nulle du fibr\'e
$\overset{2}{\bigwedge}\mathcal{T}_{X}$ d\'efinit une structure de
Poisson r\'eguli\`ere d\`es que le fibr\'e $\mathcal{F}$ coorespondant est
int\'egrable.\\}
\newline
\indent\textbf{Corollaire 1.2}\textit{ Soit $X$ une vari\'et\'e projective de dimension
  3 admettant une structure de Poisson quasi-r\'eguli\`ere.\\
\indent Alors il existe une injection de faisceaux
$\omega_{X}\hookrightarrow\Omega^{1}_{X}$ dont le conoyau est sans
torsion et localement libre de rang 2 sauf en un nombre fini de points.}\\
\newline
\indent\textbf{Proposition 1.3}\textit{ Soit $X$ une vari\'et\'e
  projective de dimension 
  3 admettant une structure de Poisson quasi-r\'eguli\`ere.\\
\indent Alors:\\
\indent(1) La premi\`ere classe de Chern de $X$, $c_{1}(X)\in
H^{1}(X,\Omega_{X}^{1})\hookrightarrow H^{2}(X,\mathbb{C})$, provient
de $H^{1}(X,\omega_{X})$ par la fl\`eche
$H^{1}(X,\omega_{X})\longrightarrow H^{1}(X,\Omega_{X}^{1})$ associ\'ee
\`a la structure de Poisson,\\
\indent(2) $c_{1}(X)^{2}= 0$ dans
$H^{4}(X,\mathbb{C})$.\\}
\newline
\indent\textit{D\'emonstration} Remarquons que le point $(2)$ est une
cons\'equence imm\'ediate du point $(1)$.
 Soit $U$ l'ouvert de Zariski de $X$ compl\'ementaire dans $X$
des points o\`u la structure de Poisson est nulle. Le ferm\'e $Y=X-U$ est
donc de codimension 3. Par suite, les groupes de cohomologie
$H_{Y}^{i}(X,\omega_{X})$ et
$H_{Y}^{i}(X,\Omega_{X}^{1})$ sont nuls pour $i\in\{0,1,2\}$ et il
en r\'esulte des isomorphismes naturels
$H^{1}(X,\omega_{X})\cong H^{1}(U,\omega_{U})$ et
$H^{1}(X,\Omega_{X}^{1})\cong H^{1}(U,\Omega_{U}^{1})$. Remarquons en
outre que la restriction \`a l'ouvert $U$ de l'injection de faisceaux
$\omega_{X}\hookrightarrow\Omega^{1}_{X}$ est une injection de fibr\'es
vectoriels dont le conoyau est naturellement un sous-fibr\'e int\'egrable de
$\mathcal{T}_{U}$ (cf. lemme 1.1). En utilisant les
arguments de R.Bott et P.Baum ([BB] prop. 3.3
et cor. 3.4), on en d\'eduit que la classe d'Atiyah du fibr\'e
$\omega_{U}$ dans $H^{1}(U,\Omega_{U}^{1})$ provient d'un \'el\'ement de
$H^{1}(U,\omega_{U})$ par la fl\`eche naturelle
$H^{1}(U,\omega_{U})\longrightarrow H^{1}(U,\Omega_{U}^{1})$. On v\'erifie
que le diagramme suivant est commutatif:
\begin{equation*}
\begin{CD}
H^{1}(X,\omega_{X})@)))H^{1}(X,\Omega_{X}^{1})\\
@VV\wr V @VV\wr V\\
H^{1}(U,\omega_{U})@)))H^{1}(U,\Omega_{U}^{1})
\end{CD}
\end{equation*}
La formation de la classe d'Atiyah \'etant fonctorielle, il en r\'esulte
que la classe d'Atiyah du fibr\'e $\omega_{X}$ dans
$H^{1}(X,\Omega_{X}^{1})$ provient d'un \'el\'ement de
$H^{1}(X,\omega_{X})$ par l'application
$H^{1}(X,\omega_{X})\longrightarrow H^{1}(X,\Omega_{X}^{1})$ ce qui
constitue la premi\`ere assertion de notre proposition, puisque la
classe d'Atiyah du fibr\'e $\omega_{X}$ et $c_{1}(X)$ sont
proportionnelles.\qed\\
\newline
\indent Rappelons qu'une vari\'et\'e est dite
minimale lorsque le fibr\'e canonique est num\'eriquement effectif.\\
\newline
\indent\textbf{Corollaire 1.4}\textit{ Soit $X$ une vari\'et\'e
  projective minimale de dimension 3 admettant une structure de Poisson
  quasi-r\'eguli\`ere.\\
\indent Alors la dimension de Kodaira $\kappa(X)$ de $X$ v\'erifie les
in\'egalit\'es $0\le\kappa(X)\le1$.}\\
\newline  
\indent\textit{D\'emonstration}  Par le th\'eor\`eme d'Abondance ([MP]), 
la dimension de Kodaira
$\kappa(X)$ d'une vari\'et\'e projective minimale de dimension 3 se calcule
num\'eriquement par la formule $\kappa(X)=\text{max}\{
i|\omega_{X}^{i}H^{3-i}\neq 0\}$ o\`u $H$ est une section hyperplane de
$X$, de sorte que, par la proposition 1.3, on obtient
les in\'egalit\'es $0\le\kappa(X)\le 1$.\qed
\vspace{1cm}\\
\centerline{{\Large\textbf{\S2 Cas o\`u $X$ n'est pas minimale}}}\\
\newline
\indent\textbf{Lemme 2.1}\textit{ Soit $X$ une vari\'et\'e projective non
  minimale admettant une structure de Poisson quasi-r\'eguli\`ere.\\
\indent Alors $X$ est un fibr\'e en coniques.}\\
\newline
\indent\textit{D\'emonstration}  Par le th\'eor\`eme de structure de S.Mori
([Mo] thm. 3.3 et thm. 3.5)
il suffit d'\'etudier les cinq cas suivants.\\
\newline
\indent\textit{Cas 1 } Il existe un diviseur $D\subset X$ tel que le
couple $(D,\mathcal{O}_{D}(D))$ soit
$(\mathbb{P}^{2},\mathcal{O}_{\mathbb{P}^{2}}(-1))$, 
$(\mathbb{P}^{2},\mathcal{O}_{\mathbb{P}^{2}}(-2))$
ou
$(\mathcal{Q},\mathcal{O}_{\mathbb{P}^{3}}(-1)\otimes\mathcal{O}_{\mathcal{Q}})$,
$\mathcal{Q}$ \'etant une quadrique r\'eduite et irr\'eductible de
$\mathbb{P}^{3}$. On en d\'eduit que ${\omega_{X}}_{|D}^{2}=4,1$ ou $2$
par la formule d'adjonction, ce qui constitue une contradiction compte
tenu de la proposition 1.3.\\
\newline
\indent\textit{Cas 2 } Le fibr\'e $\omega_{X}^{-1}$ est ample. Ce cas est
imm\'ediatement \'elimin\'e par la relation $\omega_{X}^{2}=0.$\\
\newline
\indent\textit{Cas 3 } Il existe un morphisme $X\longrightarrow Y$, $Y$
\'etant une courbe lisse, tel que pour tout point g\'eom\'etrique $\eta\in
X$, la fibre $X_{\eta}$ soit une surface int\`egre \`a fibr\'e anticanonique
$\omega_{X_{\eta}}^{-1}$ ample. Si $F$ est une fibre g\'en\'erique, on a
$\omega_{F}={\omega_{X}}_{|F}$ par la formule d'adjonction et donc $\omega_{F}^{2}=0$ ce qui contredit
l'amplitude du fibr\'e anticanonique ${\omega_{F}}^{-1}.$\\ 
\newline 
\indent\textit{Cas 4 } Il existe un morphisme $X\overset{\pi}{\longrightarrow} Y$, $Y$
\'etant une vari\'et\'e projective lisse, et une courbe lisse $C\hookrightarrow Y$
tels que $X$ soit l'\'eclat\'e de la courbe $C$ dans $Y$. L'isomorphisme
$\pi_{*}\mathcal{O}_{X}\simeq\mathcal{O}_{Y}$ permet de d\'efinir une
structure de Poisson sur $Y$ \`a l'aide de celle existant sur $X$. On
v\'erifie que le diagramme suivant est commutatif, les fl\`eches
horizontales \'etant donn\'ees par les structures de Poisson:
\begin{equation*}
\begin{CD}
H^{1}(X,\omega_{X}) @)))H^{1}(X,\Omega_{X}^{1})\\
@AA\wr A@AAA\\
H^{1}(Y,\omega_{Y}) @))) H^{1}(Y,\Omega_{Y}^{1})\\
\end{CD}
\end{equation*}
\noindent et on d\'eduit de la proposition 1.3 que $c_{1}(X)$ provient d'un
\'el\'ement de $H^{1}(Y,\Omega_{Y}^{1})$, ce qui est absurde puisque
$c_{1}(X)=-E+\pi^{*}c_{1}(\omega_{Y})$ et puisque $E$ est contract\'e
par $\pi$.\\
\newline
\indent\textit{Cas 5 } Il existe un morphisme $X\longrightarrow Y$, $Y$
\'etant une surface projective lisse et connexe, tel que, pour tout
point g\'eom\'etrique $\eta$ de $Y$, la fibre $X_{\eta}$ soit une conique
de $\mathbb{P}^{2}_{k(\eta)}.$\qed \\
\newline
\indent La fin de cette section est consacr\'ee \`a l'\'etude de ce dernier
cas.\\
\newline
\indent\textbf{Proposition 2.2}\textit{ Soit
  $X\overset{\pi}{\longrightarrow}S$ un fibr\'e en coniques admettant
  une structure de Poisson quasi-r\'eguli\`ere.\\
\indent Alors le morphisme $\pi$ est lisse et $S$ est une surface
$K3$ ou une surface ab\'elienne.\\
\indent En outre, la structure de Poisson est r\'eguli\`ere.}\\
\newline
\indent\textit{D\'emonstration} Supposons que le morphisme $\pi$ ne soit
pas lisse. Alors le lieu de d\'eg\'en\'erescence de $\pi$ est un diviseur
effectif $C_{0}$ r\'eduit, \`a croisements normaux. Au-dessus d'un point
r\'egulier de $C_{0}$, la fibre sch\'ematique de $\pi$ est la r\'eunion de
deux courbes rationnelles lisses distinctes (et donc transverses) et,
au-dessus d'un point singulier de $C_{0}$ (un point double ordinaire),
la fibre sch\'ematique de $\pi$ est une droite double ([Be2] prop. 1.2).\\
\indent Nous allons prouver la formule $4K_{S}+C_{0}\equiv0$. Soit
$C$ une courbe lisse et connexe de $S$ telle que
$C_{0}$ et $C$ se coupent tranversalement. En particulier,
les points d'intersection de $C$ et $C_{0}$ sont des points
r\'eguliers de $C_{0}$. Posons $D=\pi^{-1}(C)$. Un calcul en
coordonn\'ees locales montre que $D$ est une surface lisse (connexe). En
outre, il n'est pas difficile de voir que l'on peut contracter l'une des
composantes irr\'eductibles de chaque  fibre singuli\`ere du morphisme
$D\overset{\pi_{|D}}{\longrightarrow}C$, de sorte que $D$ est l'\'eclat\'e
en $C_{0}.C$ points d'une surface r\'egl\'ee lisse et connexe
$D_{0}$. Soit $D\overset{\phi}{\longrightarrow}D_{0}$ le morphisme
correspondant. L'application rationnelle naturelle $D_{0}\cdots\rightarrow
C$ se prolonge en un morphisme
$D_{0}\overset{\pi_{0}}{\longrightarrow}C$ tel que
$\pi_{|D}=\pi_{0}\phi$, de sorte que la surface
$D_{0}\overset{\pi_{0}}{\longrightarrow}C$ est une surface
g\'eom\'etriquement r\'egl\'ee au-dessus de $C$. La formule
d'adjonction fournit l'\'egalit\'e
${\omega_{X}}_{|D}=\omega_{D}\otimes{\mathcal{N}_{D/X}^{-1}}$. Puisque
$\mathcal{N}_{D/X}=\pi_{|D}^{*}\mathcal{N}_{C/S}$,
$\mathcal{N}_{D/X}$ est num\'eriquement \'equivalent \`a $C^{2}$
fibres du morphisme $D\overset{\pi_{|D}}{\longrightarrow}C$. On en d\'eduit que
$\omega_{D}^{2}={\omega_{{D}_{0}}}^{2}-C_{0}.C=8(1-g(C))-C_{0}.C$ et
$\omega_{D}.\mathcal{N}_{D/X}=-2C^{2}$ (cf. [H] Chap. V.2, V.3).
Par la proposition 1.3, on a ${\omega_{X}}_{|D}^{2}=0$ et on obtient
la formule $8(1-g(C))-C_{0}.C+4C^{2}=0$. Or la formule d'adjonction
donne la relation $2(g(C)-1)=C^{2}+C.K_{S}$, et on obtient finalement
la formule $(4K_{S}+C_{0}).C=0$ ce qui prouve notre assertion.\\
\indent Par la proposition 1.3, l'\'el\'ement $c_{1}(X)\in
H^{1}(X,\Omega_{X}^{1})\hookrightarrow H^{2}(X,\mathbb{C})$ est image d'un
\'el\'ement de $H^{1}(X,\omega_{X})$ par l'application
$H^{1}(X,\omega_{X})\longrightarrow H^{1}(X,\Omega_{X}^{1})$. Puisque
$c_{1}(X)$ est non nul, il en r\'esulte que l'espace vectoriel
$H^{1}(X,\omega_{X})$ est non nul. Mais, par dualit\'e de Serre, on a
$h^{1}(X,\omega_{X})=h^{2}(X,\mathcal{O}_{X})$. De plus, on v\'erifie que
$R^{1}\pi_{*}\mathcal{O}_{X}=0$ et que
$\pi_{*}\mathcal{O}_{X}=\mathcal{O}_{S}$. Il en r\'esulte en particulier
que $h^{2}(X,\mathcal{O}_{X})=h^{2}(S,\mathcal{O}_{S})$ et, par dualit\'e de
Serre \`a nouveau, on obtient finalement que $h^{0}(S,\omega_{S})\ge1$.
Le diviseur $K_{S}$ est donc effectif et l'\'egalit\'e
$4K_{S}+C_{0}\equiv0$ est impossible, ce qui d\'emontre la premi\`ere
assertion de notre proposition.\\
\indent Dans le cas o\`u le morphisme $\pi$ est lisse, la m\^eme m\'ethode
permet de prouver que $K_{S}\equiv0$ et, puisqu'on a toujours 
$h^{0}(S,\omega_{S})\ge1$, il en r\'esulte que
$\omega_{S}=\mathcal{O}_{S}$ et donc que $S$ est une surface ab\'elienne
ou une surface $K3$.\\
\indent Montrons maintenant que la structure est r\'eguli\`ere. Puisque
$\pi$ est lisse, $\pi$ est en particulier 
topologiquement localement trivial et on a donc la formule
$e(X)=e(\mathbb{P}^{1})e(S)=2e(S)$. Soit $\sigma\in
H^{0}(X,\overset{2}{\wedge}\mathcal{T}_{X})$ le bivecteur de Poisson
d\'efinissant la 
structure de Poisson sur $X$. Par hypoth\`ese, le lieu $Z$ des z\'eros de
$\sigma$ est de codimension 3 et donc
$c_{3}(\overset{2}{\wedge}\mathcal{T}_{X})=\text{deg}(Z)$. Or,
$c_{3}(\overset{2}{\wedge}\mathcal{T}_{X})=c_{1}c_{2}-c_{3}$ et on a
donc la formule $c_{1}c_{2}-c_{3}=\text{deg}(Z)$. En outre, pour un
fibr\'e en coniques, on a 
$\chi(\mathcal{O}_{X})=\chi(\mathcal{O}_{S})$, puisque
$R^{1}\pi_{*}\mathcal{O}_{X}=0$ et
$\pi_{*}\mathcal{O}_{X}=\mathcal{O}_{S}$. Enfin, par la formule
de Noether, $e(S)=12\chi(\mathcal{O}_{S})$. Finalement, on obtient
$c_{1}c_{2}-c_{3}=24\chi(\mathcal{O}_{X})-e(X)=0$ en utilisant la
formule de Riemann-Roch. Et donc $\text{deg}(Z)=0$, ce qui permet de
conclure.\qed\\
\newline
\indent\textbf{Proposition 2.3}\textit{ Soit
  $X\overset{\pi}{\longrightarrow}S$ un fibr\'e en droites projectives 
 et supposons que $X$ admette une structure de Poisson
  r\'eguli\`ere.\\
\indent  Alors la surface $S$ admet une structure de Poisson non
triviale, induite par l'isomorphisme naturel $\mathcal{O}_{S}\cong\pi_{*}\mathcal{O}_{X}$.}\\
\newline
\indent\textit{D\'emonstration} L'isomorphisme naturel
$\mathcal{O}_{S}\cong\pi_{*}\mathcal{O}_{X}$ permet de munir $S$ d'une
structure de Poisson. Supposons que cette structure soit nulle et
montrons que nous aboutissons \`a une absurdit\'e.\\
\indent On v\'erifie, par un calcul en coordonn\'ees locales, que la fl\`eche
$\omega_{X}\longrightarrow\Omega^{1}_{X/S}$ obtenue par composition
des applications $\omega_{X}\longrightarrow\Omega_{X}^{1}$ (cf. lemme
1.1) et $\Omega_{X}^{1}\longrightarrow\Omega^{1}_{X/S}$ 
est nulle et que cette assertion est \'equivalente \`a la nullit\'e de la
structure de Poisson induite sur $S$. Par suite, cette application se
factorise \`a travers
$\pi^{*}\Omega_{S}^{1}=\text{Ker}(\Omega_{X}^{1}\longrightarrow\Omega^{1}_{X/S})$
et on obtient une suite exacte:
$$0\longrightarrow\omega_{X}\longrightarrow\pi^{*}\Omega_{S}^{1}
\longrightarrow\mathcal{T}_{X/S}\longrightarrow 0.$$ 
\noindent On en d\'eduit
alors l'\'egalit\'e
$\pi^{*}c_{2}(\Omega_{S}^{1})=c_{1}(\mathcal{T}_{X/S}).c_{1}(\omega_{X})=(\pi^{*}c_{1}(\Omega_{S}^{1})-c_{1}(\omega_{X})).c_{1}(\omega_{X})$.
En cohomologie, on obtient l'\'egalit\'e
$\pi^{*}c_{2}(\Omega_{S}^{1})=\pi^{*}c_{1}(\Omega_{S}^{1}).c_{1}(\omega_{X})$
puisque $\omega_{X}^{2}\equiv 0$ dans
$H^{4}(X,\mathbb{C})$ (cf. proposition 1.3). Par la proposition
pr\'ec\'edente, nous savons aussi que $\omega_{S}\cong\mathcal{O}_{S}$ et
il en r\'esulte l'\'egalit\'e $\pi^{*}c_{2}(\Omega_{S}^{1})=0$. Rappelons que le
diviseur anticanonique $K_{X}^{-1}$ est de degr\'e relatif $-2$ au
dessus de $S$, et donc,
en utilisant la formule de projection, on en d\'eduit que
$e(S)=c_{2}(\Omega^{1}_{S})=0$ et que $S$ est une surface
ab\'elienne.\\
\indent La fl\`eche
$\pi^{*}\Omega_{S}^{1}=\mathcal{O}_{X}\oplus\mathcal{O}_{X}\longrightarrow\omega_{X}^{-1}\longrightarrow0$
fournit un syst\`eme lin\'eaire sans points bases $\Lambda\subset|
K_{X}^{-1}|$ de dimension 2 et donc un morphisme
$X\overset{f_{0}}{\longrightarrow}\mathbb{P}^{1}$ tel que
$\omega_{X}^{-1}=f_{0}^{*}\mathcal{O}_{\mathbb{P}^{1}}(1)$. Soit
$X\overset{f}{\longrightarrow}C$ la factorisation de
Stein du morphisme $f_{0}$, $C$ \'etant une courbe lisse et
connexe. Il existe donc un morphisme
$C\overset{\alpha}{\longrightarrow}\mathbb{P}^{1}$ tel que
$f_{0}=\alpha f$.
Soit $F$ une fibre g\'en\'erique lisse (connexe) de $f$. Par la
formule d'adjonction $\omega_{F}={\omega_{X}}_{|F}$ et donc
$\omega_{F}=\mathcal{O}_{F}$ puisque $\omega_{X}=f^{*}(\alpha^{*}\mathcal{O}_{\mathbb{P}^{1}}(-1))$. La formule de sous-additivit\'e de la
dimension de Kodaira ([Ka2]) entra\^{\i}ne que la dimension de Kodaira de
$C$ vaut $\kappa(C)=-\infty$ et donc que $C\cong\mathbb{P}^{1}$. Par
suite,
$\alpha^{*}\mathcal{O}_{\mathbb{P}^{1}}(1)=\mathcal{O}_{\mathbb{P}^{1}}(d)$
pour un certain entier $d\ge1$ et donc
$\omega_{X/\mathbb{P}^{1}}=f^{*}\mathcal{O}_{\mathbb{P}^{1}}(2-d)$. Par un r\'esultat
de T.Fujita ([F] prop. 1.2),
$\text{deg}_{\mathbb{P}^{1}}(\omega_{X/\mathbb{P}^{1}})\ge0$ et donc
$d\le2$. Il nous faut donc \'eliminer les deux cas $d=1$ et $d=2$.\\
\newline
\indent\textit{Cas d=2} Nous allons prouver que $X\cong
\mathbb{P}^{1}\times S$ et obtenir une contradiction. Dans ce cas, le
syst\`eme lin\'eaire $|K_{X}^{-1}|$ est de dimension 3 et
$\omega_{X}=f^{*}(\omega_{\mathbb{P}^{1}})$. Par un r\'esultat de
T.Fujita ([F], thm. 4.8), le morphisme $f$ est lisse (\`a fibres
connexes). Soit $F$ une fibre de $f$. Par la formule d'adjonction,
$\omega_{F}={\omega_{X}}_{|F}$ et donc
$\omega_{F}=\mathcal{O}_{F}$. Il en r\'esulte que $F$ est une surface
$K3$ ou une surface ab\'elienne. Mais nous savons que $S$ est
une surface ab\'elienne et il ne peut y avoir de morphisme non constant
d'une surface $K3$ vers une surface ab\'elienne. On en d\'eduit que $F$
est une surface ab\'elienne et donc qu'elle ne peut contenir aucune
fibre de $\pi$ (ce sont 
des courbes rationnelles). Par suite, le morphisme 
$F\overset{\pi_{|F}}{\longrightarrow}S$ est
fini. Si $h$ d\'esigne une fibre de $\pi$, on a $F.h=1$ et
donc la surface $F$ est une section de $\pi$. Il n'est pas difficile
de voir que le morphisme $f\times\pi$ induit un isomorphisme de $X$
sur $\mathbb{P}^{1}\times S$.\\
\indent Il nous reste \`a voir qu'il n'existe sur
$\mathbb{P}^{1}\times S$ aucune structure de Poisson induisant la
structure nulle sur $S$ via la premi\`ere projection. On a
un isomorphisme $\mathcal{T}_{X}\cong
p^{*}\mathcal{T}_{\mathbb{P}^{1}}\oplus q^{*}\mathcal{T}_{S}$ o\`u $p$
et $q$ sont respectivement les projections de $X$ sur $\mathbb{P}^{1}$
et $S$. Par suite, on a un isomorphisme
$\overset{2}{\wedge}\mathcal{T}_{X}\cong
p^{*}\mathcal{T}_{\mathbb{P}^{1}}\otimes q^{*}\mathcal{T}_{S}\oplus
q^{*}\omega_{S}^{-1}$ et donc un isomorphisme
$H^{0}(X,\overset{2}{\wedge}\mathcal{T}_{X})\cong
H^{0}(X,p^{*}\mathcal{T}_{\mathbb{P}^{1}}\otimes q^{*}\mathcal{T}_{S})
\oplus H^{0}(X,q^{*}\omega_{S}^{-1})$. La structure de Poisson sur $X$
est donn\'ee par un \'el\'ement $\sigma\in
H^{0}(X,\overset{2}{\wedge}\mathcal{T}_{X})$ qui, compte tenu des
hypoth\`eses, n'a pas de composante sur
$H^{0}(X,q^{*}\omega_{S}^{-1})$. Enfin, par la formule de K\"unneth, on
a un isomorphisme $H^{0}(X,p^{*}\mathcal{T}_{\mathbb{P}^{1}}\otimes
q^{*}\mathcal{T}_{S})\cong
H^{0}(\mathbb{P}^{1},\mathcal{T}_{\mathbb{P}^{1}})\otimes
H^{0}(S,\mathcal{T}_{S})$. Soient $(e_{1},e_{2})$ une base de
$H^{0}(S,\mathcal{T}_{S})$ et $(f_{1},f_{2},f_{3})$ une base de
$H^{0}(\mathbb{P}^{1},\mathcal{T}_{\mathbb{P}^{1}})$. On peut \'ecrire
$\sigma=\sum_{i,j}a_{i,j}e_{i}\wedge f_{j}$ o\`u les $a_{i,j}$ sont des
constantes. Choisissons des coordonn\'ees locales $(x_{1},x_{2})$ sur $S$ de
sorte que la base $(e_{1},e_{2})$ soit donn\'ee par les formules
$e_{1}=\partial_{x_{1}}$ et $e_{2}=\partial_{x_{2}}$. Enfin, on choisit une
coordonn\'ee locale $z$ sur $\mathbb{P}^{1}$ correspondant \`a une
d\'ecomposition $\mathbb{P}^{1}=\mathbb{C}\cup\{\infty\}$ de sorte que
la base $(f_{1},f_{2},f_{3})$ soit donn\'ee par les formules
$f_{1}=\partial_{z}$, $f_{2}=z\partial_{z}$ et
$f_{3}=z^{2}\partial_{z}$. Si $f$ et $g$ sont deux fonctions
holomorphes d\'efinies localement, le crochet de Poisson est alors donn\'e
par la formule
$\{f,g\}=\sum_{i,j}a_{i,j}z^{j}(\partial_{x_{i}}f\partial_{z}g-\partial_{z}f\partial_{x_{i}}g)$.
En particulier, $\{x_{1},z\}=a_{10}+a_{11}z+a_{12}z^{2}=P_{1}(z)$ et
$\{x_{2},z\}=a_{20}+a_{21}z+a_{22}z^{2}=P_{2}(z)$. Enfin, l'identit\'e
de Jacobi donne la relation $P_{1}^{'}P_{2}=P_{1}P_{2}^{'}$ et donc
les polyn\^omes $P_{1}$ et $P_{2}$ sont proportionnels, ce qui constitue
la contradiction cherch\'ee.\\
\newline
\indent\textit{Cas d=1} Dans ce cas, le syst\`eme lin\'eaire
$\Lambda=|K_{X}^{-1}|$ est de dimension 2, tous ses \'el\'ements sont
connexes et on a la relation
$\omega_{X}=f^{*}\mathcal{O}_{\mathbb{P}^{1}}(-1)$. Par le th\'eor\`eme de
lissit\'e g\'en\'erique, un \'el\'ement g\'en\'erique $D\in|K_{X}^{-1}|$ est
r\'egulier (le syst\`eme lin\'eaire est sans point base). Soit $D$ un tel
\'el\'ement. Le diviseur $D$ est irr\'eductible puisque connexe et
r\'egulier. De m\^eme que nous l'avons fait au cours de la preuve du cas
$d=2$, on v\'erifie que $D$ est une surface ab\'elienne et que le
morphisme $D\overset{\pi_{|D}}{\longrightarrow}S$ est fini. On en
d\'eduit que ce morphisme est automatiquement non ramifi\'e. On obtient
donc un rev\^etement \'etale de degr\'e
2, $D\overset{\pi_{|D}}{\longrightarrow}S$. Or, les rev\^etements \'etales
de degr\'e 2 de la surface ab\'elienne $S$ sont classifi\'es par la
2-torsion de $\text{Pic}(S)$ et sont donc en nombre fini, \`a
isomorphisme pr\`es. On peut donc trouver $D_{1}\in|K_{X}^{-1}|$ et
$D_{2}\in|K_{X}^{-1}|$ deux diviseurs lisses tels que les rev\^etements
$D_{1}\overset{\pi_{|D_{1}}}{\longrightarrow}S$ et
$D_{2}\overset{\pi_{|D_{2}}}{\longrightarrow}S$ soient
isomorphes. Effectuons le changement de base $D_{1}\longrightarrow
S$. On obtient un diagramme cart\'esien:
\begin{equation*}
\begin{CD}
\overline{X}@)p))X\\
@VVqV @VV\pi V\\
D_{1} @)\pi_{D_{1}}))S
\end{CD}
\end{equation*}
En outre, le morphisme $p$ est un rev\^etement \'etale et le morphisme $q$
poss\`ede deux sections \`a supports disjoints. Puisque $p$ est \'etale, la
vari\'et\'e $\overline{X}$ admet une structure de Poisson r\'eguli\`ere
induite par celle de $X$ et il n'est pas difficile de voir que cette
structure induit la structure nulle sur $D_{1}$ si et seulement si la
structure de Poisson sur $X$ induit la structure nulle sur $S$. On
peut donc supposer que le morphisme $\pi$ admet deux sections \`a
supports disjoints et que l'on se trouve toujours dans le cas $d=1$
puisque le cas $d=2$ est \'elimin\'e. Alors l'existence d'une section
nous assure qu'il existe un fibr\'e $\mathcal{E}$ sur $S$ localement
libre de rang 2 tel que $X=\mathbb{P}(\mathcal{E})$. En outre,
l'existence de deux sections \`a supports disjoints entra\^{\i}ne que le
fibr\'e $\mathcal{E}$ est d\'ecompos\'e, $i.e$
$\mathcal{E}=\mathcal{L}\oplus\mathcal{M}$ o\`u $\mathcal{L}$ et
$\mathcal{M}$ sont deux fibr\'es en doites. De plus, on peut toujours
supposer que $\mathcal{M}=\mathcal{O}$, de sorte que
$c_{2}(\mathcal{E})=0$. Notons $\xi$ le fibr\'e tautologique sur
$\mathbb{P}(\mathcal{E})$. Nous avons alors la formule de Grothendieck
$\xi^{2}-\pi^{*}c_{1}(\mathcal{E}).\xi+\pi^{*}c_{2}(\mathcal{E})=0$
dans $H^{4}(\mathbb{P}(\mathcal{E}),\mathbb{C})$. Et, puisque
$c_{2}(\mathcal{E})=0$, on a
$\xi^{2}=\pi^{*}c_{1}(\mathcal{E}).\xi$. On a aussi la formule
$c_{1}(X)=2\xi-\pi^{*}c_{1}(\mathcal{E})-\pi^{*}c_{1}(S)=2\xi-\pi^{*}c_{1}(\mathcal{E})$
car $c_{1}(S)=0$. Enfin, l'\'egalit\'e $c_{1}(X)^{2}=0$ entraine
$c_{1}(X)^{2}.\xi=0$. Or
$c_{1}(X)^{2}.\xi=4\xi^{3}-4\xi^{2}.\pi^{*}c_{1}(\mathcal{E})+\pi^{*}c_{1}(\mathcal{E})^{2}.\xi=c_{1}(\mathcal{E})^{2}$.
Finalement, on obtient
$c_{1}(\mathcal{E})^{2}=c_{1}(\mathcal{L})^{2}=0$.\\
\indent Nous allons prouver qu'il existe une courbe elliptique
$E\subset S$ telle que $\mathcal{L}\cong\mathcal{O}(E)$. En effet, on
a un isomorphisme
$\omega_{X}^{-1}\cong\xi^{2}\otimes\pi^{*}(\text{det}(\mathcal{E})^{-1})$
et donc
$H^{0}(X,\omega_{X}^{-1})=H^{0}(S,S^{2}(\mathcal{E})\otimes\text{det}(\mathcal{E})^{-1}$.
Or,
$S^{2}(\mathcal{E})\otimes\text{det}(\mathcal{E})^{-1}=\mathcal{O}_{S}\oplus\mathcal{L}\oplus\mathcal{L}^{-1}$
et $h^{0}(X,\omega_{X}^{-1})=2$. Il en r\'esulte que
$h^{0}(S,\mathcal{L})+h^{0}(S,\mathcal{L}^{-1})=1$. Supposons par
exemple que $h^{0}(S,\mathcal{L})=1$. Alors,
$\mathcal{L}=\mathcal{O}_{S}(D)$ o\`u $D$ est un diviseur effectif de
$S$. Ecrivons $D=\sum_{i}n_{i}C_{i}$ o\`u les $C_{i}$ sont les
composantes irr\'eductibles de $D$. Comme $S$ est une surface ab\'elienne,
$C_{i}^{2}\ge0$. Comme $D^{2}=0$, on en d\'eduit donc que
$C_{i}.C_{j}=0$ pour tout couple d'entier $(i,j)$. Il r\'esulte alors de
la formule d'adjonction que $p_{a}(C_{i})=1$ et, comme il n'y a pas de
morphisme non constant $\mathbb{P}^{1}\longrightarrow S$, on en d\'eduit
que $C_{i}$ est une courbe elliptique. Puisque $C_{i}.C_{j}=0$, les
$C_{i}$ sont disjoints et puisque $h^{0}(S,\mathcal{O}(D))=1$,
$D$ est un diviseur r\'eduit et irr\'eductible et donc une
courbe elliptique que nous noterons $E$.\\
\indent Par un r\'esultat de Poincar\'e
([Mum] p.173) il existe une courbe elliptique $B\subset S$ telle
que $B\cap E$ soit fini. Par suite, le morphisme $B\times
E\overset{(x,y)\mapsto x+y}{\longrightarrow}
S$ est un rev\^etement \'etale fini et on peut toujours supposer que ce
morphisme est un morphisme de groupes. Effectuons alors le changement
de base $B\times E\longrightarrow S$, de sorte que l'on obtient le
diagramme cart\'esien suivant: 
\begin{equation*}
\begin{CD}
\overline{X}@)p))X\\
@VVqV @VV\pi V\\
B\times E @)\alpha))S
\end{CD}
\end{equation*}
Comme $p$ est un rev\^etement \'etale, la vari\'et\'e $\overline{X}$ admet une
structure de Poisson r\'eguli\`ere induite par celle de $X$ et il n'est
pas difficile de voir que cette structure induit la structure nulle
sur $B\times E$ si et seulement si la structure de Poisson sur $X$
induit la structure nulle sur $S$. Le cas $d=2$ \'etant \'elimin\'e, nous
pouvons supposer que l'on se trouve dans le cas $d=1$. Or
$\overline{X}$ est le fibr\'e projectif associ\'e au fibr\'e
$\alpha^{*}\mathcal{E}=\mathcal{O}_{B\times
  E}\oplus\mathcal{O}_{B\times E}(\alpha^{-1}(E))$ et $\alpha^{-1}(E)$
est la r\'eunion disjointe des $\{x\}\times E$ pour $x\in B\cap E$. Par
les calculs d\'eja faits on a $h^{0}(B\times E,\mathcal{O}_{B\times
  E}(\alpha^{-1}(E)))+h^{0}(B\times E,\mathcal{O}_{B\times
  E}(-\alpha^{-1}(E)))=1$ et donc $h^{0}(B\times E,\mathcal{O}_{B\times
  E}(\alpha^{-1}(E)))=1$, ce qui entraîne que $\alpha^{-1}(E)$ est
r\'eduit \`a $0_{B}\times E$ et que $B\cap E=0_{S}$. Aussi, $\alpha$ est
en fait un isomorphisme et il nous reste donc \`a traiter le cas o\`u
$X\cong\mathbb{P}_{E\times B}(\mathcal{O}_{E\times
  B}\oplus\mathcal{O}_{E\times B}(E\times \{0_{B}\}))\cong E\times
\mathbb{P}_{B}(\mathcal{O}_{B}\oplus\mathcal{O}_{B}(0_{B}))\cong
E\times \mathbb{P}_{B}(\mathcal{O}_{B}\oplus\mathcal{O}_{B}(-0_{B}))$.
Remarquons que la surface
$\mathbb{P}_{B}(\mathcal{O}_{B}\oplus\mathcal{O}_{B}(-0_{B}))$ est
normalis\'ee au sens de [H] (Chap. V.2) de sorte que le fibr\'e
tautologique
$\mathcal{O}_{\mathbb{P}_{B}(\mathcal{O}_{B}\oplus\mathcal{O}_{B}(-0_{B}))}(1)$
est effectif. Il en r\'esulte que le fibr\'e anticanonique de
la surface
$\mathbb{P}_{B}(\mathcal{O}_{B}\oplus\mathcal{O}_{B}(-0_{B}))$ est
effectif et donn\'e par la formule
$K_{\mathbb{P}_{B}(\mathcal{O}_{B}\oplus\mathcal{O}_{B}(-0_{B}))}^{-1}=\mathcal{O}_{\mathbb{P}_{B}(\mathcal{O}_{B}\oplus\mathcal{O}_{B}(-0_{B}))}(2)+\mathbb{P}^{1}\times\{0_{B}\}$
de sorte qu'il existe dans le syst\`eme lin\'eaire $|K_{X}^{-1}|$ un
\'el\'ement ayant une composante verticale relativement au morphisme
$X=E\times
\mathbb{P}_{B}(\mathcal{O}_{B}\oplus\mathcal{O}_{B}(-0_{B}))\longrightarrow
S=E\times B$. La contradiction cherch\'ee est l\`a puisque le syst\`eme
lin\'eaire $|K_{X}^{-1}|$ est suppos\'e sans points bases et de dimension
2 ($h^{0}(X,K_{X}^{-1})=2$) et donc ses \'el\'ements sont disjoints
(n'oublions pas que le syst\`eme lin\'eaire $|K_{X}^{-1}|$ est de degr\'e
relatif 2 par rapport \`a $S$).\qed\\
\newline
\indent\textbf{Proposition 2.4}\textit{ Soit
  $X\overset{\pi}{\longrightarrow}S$ un fibr\'e en droites projectives. On
  suppose que $X$ admet une structure de Poisson r\'eguli\`ere induisant une
  structure de Poisson non nulle sur $S$ via l'isomorphisme
  $\mathcal{O}_{S}\cong\pi_{*}\mathcal{O}_{X}$.\\
\indent Alors $X=S\times\mathbb{P}^{1}$ o\`u $S$ est une surface $K3$,
ou bien $X$ est un fibr\'e plat en droites projectives sur $S$ une surface
ab\'elienne.}\\
\newline
\indent\textit{D\'emonstration} Nous savons d\'eja que $S$ est une surface
$K3$ ou une surface ab\'elienne par la proposition 2.2. Par suite, la
structure de Poisson induite sur $S$ est partout non nulle. Un calcul
en coordonn\'ees locales permet de se rendre compte que la fl\`eche
$\mathcal{T}_{X/S}\longrightarrow \omega_{X}^{-1}$ obtenue par
composition des fl\`eches
$\mathcal{T}_{X/S}\longrightarrow\mathcal{T}_{X}\longrightarrow
\omega_{X}^{-1}$ (cf. lemme 1.1) est un isomorphisme. 
On en d\'eduit l'existence d'une section
$\pi^{*}\mathcal{T}_{S}\overset{s}{\longrightarrow}\mathcal{T}_{X}$ de
l'application naturelle $\mathcal{T}_{X}\longrightarrow\pi^{*}\mathcal{T}_{S}$,
telle que $s(\pi^{*}\mathcal{T}_{S})=\mathcal{F}$, les notations
\'etant celles du lemme 1.1. La distribution d\'efinie par cette
section est donc int\'egrable par le lemme 1.1. 
Utilisant alors le th\'eor\`eme de
Frobenius et la compacit\'e des fibres de $\pi$, on en d\'eduit qu'il
existe un recouvrement $(U_{\alpha})_{\alpha}$ de $S$ et des isomorphismes
$\phi_{\alpha}:\,\pi^{-1}(U_{\alpha})\cong
U_{\alpha}\times\mathbb{P}^{1}$ au dessus de $U_{\alpha}$ qui
transforment la structure de Poisson sur $X$ en la structure de
Poisson sur $U_{\alpha}\times\mathbb{P}^{1}$ produit de la structure
symplectique sur $U_{\alpha}$ et de la structure nulle sur
$\mathbb{P}^{1}$. Fixons un point $x\in X$. Soit
$z_{\alpha,1},z_{\alpha,2},z_{\alpha,3}$ des coordonn\'ees locales au voisinage
de $\phi_{\alpha}(x)\in U_{\alpha}\times\mathbb{P}^{1}$ centr\'ees en
$\phi_{\alpha}(x)$. Ecrivons
$\phi_{\alpha\beta}(z_{\alpha,1},z_{\alpha,2},z_{\alpha,3})=(z_{\beta,1},z_{\beta,2},\varphi_{\alpha\beta}(z_{\alpha,1},z_{\alpha,2},z_{\alpha,3}))$.
Les fonctions de transition \'etant des morphismes de Poisson, on
obtient les \'egalit\'es:
$$
\left\lbrace
\begin{array}{l}
\phi_{\alpha\beta}^{*}\{z_{\beta,1},z_{\beta,3}\}=\{\phi_{\alpha\beta}^{*}z_{\beta,1},\phi_{\alpha\beta}^{*}z_{\beta,3}\}=\{z_{\alpha,1},\varphi_{\alpha\beta}(z_{\alpha,1},z_{\alpha,2},z_{\alpha,3})\}=\{z_{\alpha,1},z_{\alpha,2}\}\frac{\partial\varphi_{\alpha\beta}}{\partial
  z_{\alpha,2}}
  \\
\phi_{\alpha\beta}^{*}\{z_{\beta,2},z_{\beta,3}\}=\{\phi_{\alpha\beta}^{*}z_{\beta,2},\phi_{\alpha\beta}^{*}z_{\beta,3}\}=\{z_{\alpha,2},\varphi_{\alpha\beta}(z_{\alpha,1},z_{\alpha,2},z_{\alpha,3})\}=\{z_{\alpha,2},z_{\alpha,1}\}\frac{\partial\varphi_{\alpha\beta}}{\partial
  z_{\alpha,1}}
\end{array}
\right.
$$
\noindent o\`u $\{.,.\}$ d\'esigne le crochet de Poisson. Puisque la
structure de Poisson sur $S$ est partout non nulle, il en r\'esulte
finalement les deux \'equations:
$$
\left\lbrace
\begin{array}{l}
\frac{\partial\varphi_{\alpha\beta}}{\partial
  z_{\alpha,2}}\equiv0\\ 
\frac{\partial\varphi_{\alpha\beta}}{\partial z_{\alpha,1}}\equiv0,
\end{array}
\right.
$$
\noindent et les fonctions de transitions du fibr\'e
  $X\overset{\pi}{\longrightarrow}S$ sont donc localement
  constantes. Rappelons enfin qu'une surface $K3$ est simplement
  connexe et donc qu'un tel
  fibr\'e est alors trivial, ce qui ach\`eve
  la preuve de notre proposition.\qed\\
\newline
\indent\textbf{Proposition 2.5}\textit{ Soit $S$ une surface $K3$ et $A$
  une surface ab\'elienne.\\
\indent Alors $S\times\mathbb{P}^{1}$ admet une structure de Poisson
r\'eguli\`ere et il en est de m\^eme pour tout fibr\'e plat en droites
projectives sur $A$.}\\
\newline
\indent\textit{D\'emonstration} La premi\`ere assertion est claire. Soit
$X\overset{\pi}{\longrightarrow}S$ un fibr\'e plat en droites projectives sur
$A$. Soit $(U_{\alpha})_{\alpha}$ un recouvrement ouvert de $S$ et
$\phi_{\alpha}$ des isomorphismes $\pi^{-1}(U_{\alpha})\cong
U_{\alpha}\times\mathbb{P}^{1}$ au dessus de $U_{\alpha}$ de sorte que
les fonctions de transition $\phi_{\alpha}\phi_{\beta}^{-1}$ soient
localement constantes. On munit le 
produit $U_{\alpha}\times\mathbb{P}^{1}$ de la structure de Poisson
produit, une structure symplectique \'etant fix\'ee au pr\'ealable sur
$S$. On v\'erifie sans peine que les isomorphismes $\phi_{\alpha}$ sont
des morphismes de Poisson, les calculs \'etant analogues \`a ceux de la
  proposition pr\'ec\'edente, de sorte que l'on obtient sur $X$ une
structure de Poisson r\'eguli\`ere, ce qui termine la preuve de notre
proposition.\qed\\
\newline
\indent Nous avons donc d\'emontr\'e le th\'eor\`eme:\\
\newline
\indent\textbf{Th\'eor\`eme 2.6}\textit{ Soit $X$ une vari\'et\'e projective
de dimension 3, non minimale.\\
\indent Alors $X$ admet une structure de Poisson quasi-r\'eguli\`ere si et
seulement si $X=S\times\mathbb{P}^{1}$ o\`u $S$ est une surface $K3$,
ou bien $X$ est un fibr\'e en droites projectives sur $A$ surface ab\'elienne.\\
\indent En outre, une telle structure est r\'eguli\`ere.}
\vspace{1cm}\\
\centerline{{\Large\textbf{\S3 Cas o\`u $X$ est minimale et $\kappa(X)=0$}}}\\
\newline
\indent La formule num\'erique donnant
la dimension de Kodaira permet d'affirmer que $c_{1}(X)\equiv 0$ et
donc que le fibr\'e canonique $\omega_{X}$ est de torsion puisque $X$
est minimale.\\
\newline
\indent\textbf{Proposition 3.1}\textit{ Soit $X$ une vari\'et\'e projective
  telle que le fibr\'e canonique $\omega_{X}$
  soit de torsion.\\
\indent Alors les champs de vecteurs non nuls sur $X$ ne
sont pas tangents \`a la fibration d'Albanese de $X$.\\
\indent En outre, on a $q(X)\in\{0,1,2,3\}$ et, si
$q(X)=3$, alors $X$ est une vari\'et\'e ab\'elienne.}\\
\newline
\indent\textit{D\'emonstration} Par un th\'eor\`eme de Kawamata ([Ka1] thm. 8.3), le
morphisme d'Albanese $X\overset{\alpha_{X}}{\longrightarrow}\text{Alb}(X)$
est surjectif, \`a fibres lisses et connexes. En outre, il existe un
rev\^etement \'etale fini $B\overset{r}{\longrightarrow}\text{Alb}(X)$ tel
  que le produit fibr\'e 
$B{\times}_{\text{Alb}(X)}X$ soit isomorphe au dessus de $B$ au produit
$B\times F$, $F$ \'etant une fibre du morphisme d'Albanese. De plus, on
peut toujours supposer que le rev\^etement \'etale
$B\overset{r}{\longrightarrow}\text{Alb}(X)$ est
galoisien, de groupe $G$, de sorte que $X$ est isomorphe au quotient
$B\times F/G$, $G$ op\'erant sur $B\times F$ en respectant l'action de
$G$ sur $B$. Remarquons que, puisque $B/G\cong\text{Alb}(X)$, le groupe $G$
est form\'e de translations de $B$. Notons $p$ et $q$ les projections de
  $B\times F$ sur $X$ et $B$ respectivement. Consid\'erons la suite
  exacte de cohomologie:
$$0\longrightarrow H^{0}(X,\mathcal{T}_{X/\text{Alb}(X)})\longrightarrow
H^{0}(X,\mathcal{T}_{X})\longrightarrow
H^{0}(X,\alpha_{X}^{*}\mathcal{T}_{\text{Alb}(X)})\longrightarrow\cdots.$$
Il suffit donc
de prouver que l'application $H^{0}(X,\mathcal{T}_{X})\longrightarrow
H^{0}(X,\alpha_{X}^{*}\mathcal{T}_{\text{Alb}(X)})$ est surjective, puisque
ces deux espaces vectoriels ont m\^eme dimension ([Ka1] cor. 8.6).\\ 
\indent Puisque le morphisme $B\times F\overset{p}{\longrightarrow}X$
est un rev\^etement \'etale fini, il existe une application injective
$H^{0}(X,\mathcal{T}_{X})\hookrightarrow H^{0}(B\times
F,p^{*}\mathcal{T}_{X})=H^{0}(B\times F,\mathcal{T}_{B\times F})$ que
nous noterons $p^{*}$. On v\'erifie alors sans difficult\'es que le
diagramme suivant est commutatif:
\begin{equation*}
\begin{CD}
H^{0}(X,\mathcal{T}_{X}) @){p^{*}})) H^{0}(B,\mathcal{T}_{B})\oplus
H^{0}(F,\mathcal{T}_{F})\\ 
@VV{{\alpha_{X}}_{*}}V @VVV\\
H^{0}(\text{Alb}(X),\mathcal{T}_{\text{Alb}(X)}) @){\sim})) H^{0}(B,\mathcal{T}_{B})
\end{CD}
\end{equation*}
\indent Soit donc $\theta_{0}\in
H^{0}(\text{Alb}(X),\mathcal{T}_{\text{Alb}(X)})$
et posons
 $\theta_{1}=\frac{1}{|G|}\sum_{g\in G}g_{*}\theta_{0}\in
H^{0}(B\times F,\mathcal{T}_{B\times
  F})^{G}=H^{0}(X,\mathcal{T}_{X})$, o\`u nous avons identifi\'e les
espaces $H^{0}(\text{Alb}(X),\mathcal{T}_{\text{Alb}(X)})$ et
$H^{0}(B,\mathcal{T}_{B})$. Le groupe $G$ agissant de mani\`ere
compatible sur $B$ et $B\times F$, l'image de 
$\theta_{1}$ par l'application $H^{0}(X,\mathcal{T}_{X})\longrightarrow
H^{0}(\text{Alb}(X),\mathcal{T}_{\text{Alb(X)}})$ est $\theta_{0}$, ce
qui permet de conclure.\\
\indent La surjectivit\'e du morphisme d'Albanese entra\^{\i}ne
$q(X)\in\{0,1,2,3\}$ et, si
$q(X)=3$, alors $\alpha_{X}$ est un rev\^etement \'etale
\`a fibres connexes et donc un isomorphisme, ce qui ach\`eve la
d\'emonstration de notre proposition.\qed\\
\newline
\indent\textbf{Proposition 3.2}\textit{ Soit $X$ une vari\'et\'e
  projective minimale de dimension 3 et de
  dimension de Kodaira $\kappa(X)=0$. Supposons en
  outre que $X$ admette une structure de Poisson quasi-r\'eguli\`ere.\\
\indent Alors $q(X)\neq0$ et la structure est r\'eguli\`ere. En outre, si
$q(X)=1$, le fibr\'e canonique est trivial.}\\
\newline
\indent\textit{D\'emonstration }  D\'emontrons la premi\`ere assertion par
l'absurde et supposons donc que $q(X)=0$. Par le th\'eor\`eme de
Riemann-Roch,
$\chi(\mathcal{O}_{X})=\frac{1}{24}c_{1}(X)c_{2}(X)=0$. L'hypoth\`ese
entraine 
$h^{0}(X,\Omega^{1}_{X})=0$ par th\'eorie de Hodge. Compte tenu de
l'injection $\omega_{X}\hookrightarrow\Omega_{X}^{1}$, on en d\'eduit
que $h^{0}(X,\omega_{X})=0$ et donc que $h^{3}(X,\mathcal{O}_{X})=0$ par
dualit\'e de Serre. L'\'egalit\'e $\chi(\mathcal{O}_{X})=0$ s'\'ecrit alors
$1+h^{2}(X,\mathcal{O}_{X})=0$, ce qui constitue la contradiction
cherch\'ee.\\
\indent Montrons maintenant que la structure est r\'eguli\`ere. Soit $\sigma\in
H^{0}(X,\overset{2}{\wedge}\mathcal{T}_{X})$ la section d\'efinissant la
structure de Poisson sur $X$. Par hypoth\`ese, le lieu $Z$ des z\'eros de
$\sigma$ est de codimension 3 et donc
$c_{3}(\overset{2}{\wedge}\mathcal{T}_{X})=\text{deg}(Z)$. Or,
$c_{3}(\overset{2}{\wedge}\mathcal{T}_{X})=c_{1}c_{2}-c_{3}$ et on a
donc la formule $\text{deg}(Z)=c_{1}c_{2}-c_{3}=-c_{3}=-e(X)$. 
 Mais, puisque le morphisme d'Albanese
 $X\overset{\alpha_{X}}{\longrightarrow}\text{Alb}(X)$ 
est lisse, il est particulier topologiquement localement
trivial et donc $e(X)=0$ puisque $q(X)\neq0$. Il en r\'esulte que  
$\text{deg}(Z)=0$, ce qui permet de conclure.\\
\indent Supposons enfin que $q(X)=1$. Par la proposition 3.1,
$h^{0}(X,\mathcal{T}_{X/\text{Alb}(X)})=0$ et la suite exacte:
$$0\longrightarrow\wedge^{2}\mathcal{T}_{X/\text{Alb}(X)}=\omega_{X}^{-1}\longrightarrow\wedge^{2}\mathcal{T}_{X}\longrightarrow\mathcal{T}_{X/\text{Alb}(X)}\longrightarrow0$$
\noindent fournit donc l'\'egalit\'e
$h^{0}(X,\omega_{X}^{-1})=h^{0}(X,\wedge^{2}\mathcal{T}_{X})$. Par
hypoth\`ese, $h^{0}(X,\wedge^{2}\mathcal{T}_{X})\ge1$ et on en d\'eduit
l'in\'egalit\'e 
$h^{0}(X,\omega_{X}^{-1})\ge1$. Mais nous savons aussi que
$c_{1}(X)\equiv0$ et on en d\'eduit donc que
$\omega_{X}=\mathcal{O}_{X}$, ce qui constitue la derni\`ere assertion
de notre proposition.\qed\\ 
\newline
\indent\textbf{Proposition 3.3}\textit{ Soit $X$ une vari\'et\'e projective de
  dimension 3, \`a fibr\'e
  canonique trivial et v\'erifiant $q(X)=1$.\\
\indent Alors $X$
appartient \`a l'une des deux familles suivantes:\\
\indent (1) $X=(C\times S)/G$ o\`u $C$ est une courbe elliptique,
 $S$ est une surface $K3$ et $G$ un groupe fini 
 de translations de $C$
 op\'erant sur $S$ en respectant la structure symplectique,\\
\indent (2) $X=(C\times A)/G$ o\`u $C$ est une courbe elliptique,
 $A$ est une surface ab\'elienne  
et $G$ un groupe fini de translations de $C$ op\'erant sur $C\times A$
par la formule: 
$$g.(c,a)=(g.c,t_{g}(c,a)+u_{g}(a)),\,g\in G,\,c\in C,\,a\in A,$$
o\`u $u_{g}$ est
un automorphisme de groupes respectant la structure
symplectique de $A$ et $t_{g}$ une fonction r\'eguli\`ere sur $C\times A$ \`a
valeurs dans $A$.\\
\indent De plus, toute vari\'et\'e de l'une des deux familles admet une
structure de Poisson r\'eguli\`ere.}\\
\newline
\indent\textit{D\'emonstration } Le
morphisme d'Albanese $X\overset{\alpha_{X}}{\longrightarrow}\text{Alb}(X)$
\'etant surjectif, lisse, et \`a fibres connexes, on a la suite exacte
courte de fibr\'es vectoriels:
$$0\longrightarrow{{(\alpha_{X})}^{*}\omega_{\text{Alb}(X)}}=\mathcal{O}_{X}\longrightarrow\Omega_{X}^{1}\longrightarrow\Omega_{X/\text{Alb}(X)}^{1}\longrightarrow0.$$
\noindent L'injection de fibr\'es vectoriels
$\omega_{X}\cong\mathcal{O}_{X}\hookrightarrow\Omega_{X}^{1}$ fournit une
section globale du fibr\'e vectoriel
$\overset{2}{\wedge}\mathcal{T}_{X}\cong\Omega_{X}^{1}$. Cette section
d\'efinira une structure de Poisson r\'eguli\`ere sur $X$ \`a condition que
l'identit\'e de Jacobi soit v\'erifi\'ee. Le conoyau de l'injection 
$\omega_{X}\hookrightarrow\Omega_{X}^{1}$ est le fibr\'e
$\Omega_{X/\text{Alb}(X)}^{1}$, et par le lemme 1.1, il suffit donc de
prouver que le fibr\'e
$\mathcal{T}_{X/\text{Alb}(X)}\hookrightarrow\mathcal{T}_{X}$
est int\'egrable, ce qui est imm\'ediat puisque le morphisme d'Albanese
est lisse.\\
\indent D\'emontrons maintenant la seconde assertion de notre
proposition. Nous savons qu'il existe un rev\^etement \'etale galoisien de
groupe $G$, $C\longrightarrow\text{Alb}(X)$ tel que la fibration 
$C\times_{\text{Alb}(X)}X/C$ soit isomorphe \`a la fibration $C\times A/C$ o\`u
$A$ est soit une vari\'et\'e ab\'elienne, soit une surface $K3$. 
 Le groupe $G$ est \'evidemment un sous-groupe du groupe des
automorphismes de $C$ et il imm\'ediat que $G$ op\`ere sur $C$ par
translation. Remarquons alors que la structure de Poisson que nous
venons de construire sur $X$ se rel\`eve \`a $C\times A$ en une
structure de Poisson r\'eguli\`ere qui est le produit de la structure
nulle sur $C$ et d'une structure symplectique sur $A$. Pour s'en
convaincre, il suffit d'examiner le feuilletage d\'efini par cette
structure. On constate ensuite que le groupe $G$ respecte la structure
de Poisson sur $C\times A$. Enfin, la structure de l'action de $G$
r\'esulte du fait que le groupe des automorphismes d'une surface $K3$
ainsi que le groupe des automorphismes de groupes d'une surface
ab\'elienne sont discrets. Un calcul \'el\'ementaire permet de v\'erifier que
la stucture symplectique sur $A$ est respect\'ee, ce qui ach\`eve la preuve
de notre proposition, puisque la derni\`ere assertion est \'evidente.\qed\\
\newline
\indent\textbf{Proposition 3.4}\textit{ Soit $X$ une vari\'et\'e
  projective minimale de dimension 3 et de 
  dimension de Kodaira $\kappa(X)=0$. On suppose que
  $q(X)=2$.\\
\indent Alors $X$ est de
la forme $(F\times B)/G$ o\`u $F$ est une courbe elliptique, $B$ est
une surface ab\'elienne et $G$ un groupe fini de translations de $B$
op\'erant sur $F$ de sorte que
$F/G=\mathbb{P}^{1}$ et sur $F\times B$ par le produit de ses actions
sur chacun des facteurs.\\
\indent De plus, toute vari\'et\'e de cette forme admet une structure de
Poisson r\'eguli\`ere.}\\ 
\newline
\indent\textit{D\'emonstration } Nous savons que le
morphisme d'Albanese $X\overset{\alpha_{X}}{\longrightarrow}\text{Alb}(X)$
est surjectif, \`a fibres lisses et connexes et qu'il existe un
rev\^etement \'etale fini $B\longrightarrow\text{Alb}(X)$ galoisien de
groupe $G$ tel que le produit fibr\'e 
$X{\times}_{\text{Alb}(X)}B$ soit isomorphe au-dessus de $B$ au produit
$F\times B$, $F$ \'etant une fibre du morphisme d'Albanese. Il en
r\'esulte que $X$ est isomorphe au quotient
$B\times F/G$, $G$ op\'erant sur $B\times F$ en respectant l'action de
$G$ sur $B$. Remarquons enfin que, puisque $B/G\cong\text{Alb}(X)$, le
groupe $G$ 
est form\'e de translations de $B$ et qu'il est \'evident
que $G$ est un sous-groupe du groupe des automorphismes de
$B$. Utilisant un lemme de Beauville
([Be1] lemme VI 10, l'argument donn\'e dans le cas des surfaces s'adapte
facilement \`a notre situation),
 on peut toujours supposer que le groupe $G$ op\`ere
sur la courbe elliptique $F$ et que l'action de $G$ sur le produit
$F\times B$ est le produit des actions de $G$ sur $F$ et $B$. Enfin,
si $F/G$ est une courbe elliptique, alors $G$ op\`ere sur
$F$ par translations et il en r\'esulte que $X$ est une vari\'et\'e
ab\'elienne, ce qui impossible puisque $q(X)=2$.\\
\indent La derni\`ere assertion \'etant imm\'ediate, notre proposition est
d\'emontr\'ee.\qed\\
\newline
\indent Nous avons donc prouv\'e le th\'eor\`eme suivant:\\
\newline
\indent\textbf{Th\'eor\`eme 3.5}\textit{ Soit $X$ une vari\'et\'e projective
  de dimension 3. On suppose en outre que $X$ est de dimension de 
  Kodaira $\kappa(X)=0$.\\
\indent Alors $X$ admet une structure de Poisson r\'eguli\`ere si et
seulement si $X$ est minimale et v\'erifie l'une des trois conditions suivantes:\\
\indent (1) $q(X)=3$ (auquel cas $X$ est un vari\'et\'e ab\'elienne)\\
\indent (2) $q(X)=2$\\
\indent (3) $q(X)=1$ et $p_{g}(X)=1$.}
\vspace{1cm}\\
\centerline{{\Large\textbf{\S4 Cas o\`u $X$ est minimale et $\kappa(X)=1$}}}\\
\newline
\indent\textbf{D\'efinition 4.1}\textit{ Soient $X$ une vari\'et\'e alg\'ebrique,
   $C$ une courbe alg\'ebrique et
  $X\overset{f}{\longrightarrow}C$ un morphisme. On d\'efinit le diviseur de
  ramification de $f$ par la formule:
$$D_{f}=\sum_{P\in C}f^{*}P-{(f^{*}P)}_{red},$$
qui a un sens par le th\'eor\`eme de lissit\'e g\'en\'erique.}\\
\newline
\indent La proposition suivante est essentiellement due \`a M. Reid
([R]).\\
\newline
\indent\textbf{Proposition 4.2}\textit{ Soit $X$ une vari\'et\'e projective
  et soit $L\subset \Omega_{X}^{1}$ un faisceau
  inversible tel que le faisceau quotient $\Omega_{X}^{1}/L$ soit sans
  torsion. On suppose qu'il existe un entier $n\ge2$ et un syst\`eme
  lin\'eaire $\Lambda\subset |L^{n}|$ sans points bases tel que
  $\phi(X)$ soit une courbe, o\`u
  $X\overset{\phi}{\longrightarrow}\mathbb{P}^{N}$ est le morphisme
  associ\'e au syst\`eme lin\'eaire $\Lambda$.\\ 
\indent Alors il existe un morphisme surjectif, \`a fibres connexes
  $X\overset{f}{\longrightarrow}C$ vers une courbe lisse, tel que
  $f^{*}\omega_{C}\subset L\subset \Omega_{X}^{1}$, o\`u l'inclusion
  $f^{*}\omega_{C}\subset \Omega_{X}^{1}$ est l'inclusion
  naturelle. En outre, on a l'\'egalit\'e
  $L=f^{*}\omega_{C}(D_{f})$, $D_{f}$ \'etant le diviseur de
  ramification de $f$, et $f_{*}L=\omega_{C}$.}\\
\newline
\indent Rappelons le lemme suivant qui nous sera utile par la
suite.\\
\newline
\indent\textbf{Lemme 4.3 ([Be1] lemme VI 7 bis)}\textit{ Soient
  $\triangle\subset\mathbb{C}$ le 
  disque unit\'e, $U$ une vari\'et\'e analytique (non compacte) lisse de
  dimension 3 et
  $p\,:\,U\longrightarrow \triangle$ un morphisme. On suppose que
  $p^{*}0=nD$ o\`u $n\ge1$ et $D$ est un diviseur effectif r\'eduit. Soit
  $\triangle\overset{q}{\longrightarrow}\triangle$ le morphisme
  $z\longmapsto z^{n}$, $\tilde{U}$ le produit fibr\'e d\'efini par le
  diagramme:
\begin{equation*}
\begin{CD}
\tilde{U} @)\tilde{q}))U \\
 @VV\tilde{p}V @VVpV\\
\triangle @)q))\triangle
\end{CD}
\end{equation*}
$\overline{U}$ la normalisation de $\tilde{U}$, $\overline{p}$,
$\overline{q}$ les projections de $\overline{U}$ sur $\triangle$ et
$U$.\\
\indent Alors, $\overline{U}$ est lisse, $\overline{q}$ est \'etale fini
et la fibre ${\overline{p}}^{*}0$ est r\'eduite.}\\
\newline
\indent\textbf{Proposition 4.4}\textit{ Soit $X$ une vari\'et\'e projective
  de dimension 3, admettant une structure de Poisson
  quasi-r\'eguli\`ere. On suppose que $X$ est de dimension de Kodaira
  $\kappa(X)=1$.\\ 
\indent Sous ces hypoth\`eses, il existe un morphisme surjectif, \`a
fibres connexes,
$X\overset{f}{\longrightarrow}C$ vers une courbe alg\'ebrique $C$, dont les
fibres lisses sont soit des surfaces $K3$, soit des surfaces
ab\'eliennes.\\
\indent Si l'une des fibres est une surfaces $K3$, le morphisme
$f$ est lisse, analytiquement localement trivial et on a la formule
$\omega_{X}=f^{*}\omega_{C}$.\\
\indent Si l'une des fibres est une surface ab\'elienne, les
fibres lisses sont des surfaces ab\'eliennes et les fibres non lisses 
sont multiples de surfaces ab\'eliennes et on a l'\'egalit\'e 
$\omega_{X}=f^{*}\omega_{C}(D_{f})$, $D_{f}$ \'etant le
diviseur de ramification de $f$.\\
\indent Enfin, la structure de Poisson est r\'eguli\`ere.}\\
\newline
\indent\textit{D\'emonstration}  L'existence d'un morphisme
$X\overset{f}{\longrightarrow}C$ surjectif, \`a fibres connexes, vers
une courbe lisse tel que $f^{*}\omega_{C}(D_{f})=\omega_{X}$, $D_{f}$
\'etant le diviseur de ramification de $f$, r\'esulte du corollaire 1.2, de
la proposition 4.2 et du th\'eor\`eme principal de Y.Kawamata ([Ka3]).\\
\indent Ecrivons $D_{f}=D_{P_{1}}+\cdots+D_{P_{l}}$ ($l\ge1$) o\`u
$D_{P_{i}}=f^{*}P_{i}-{(f^{*}P_{i})}_{red}$, les points $P_{i}\in C$
\'etant tous distincts. Soit $f^{*}P_{i}=\sum_{j}n_{i,j}D_{i,j}$ la
d\'ecomposition de la fibre sch\'ematique $f^{*}P_{i}$ en somme de ses
composantes irr\'eductibles ($n_{i,j}\ge1$). Par le th\'eor\`eme de Bertini,
on peut trouver une section hyperplane $H$ de $X$ telle que:\\
\indent (1) $H$ soit lisse et connexe,\\
\indent (2) $H\nsupseteq D_{i,j}\cap D_{i,k}$ pour tout triplet
$(i,j,k)$,\\
\indent (3) le morphisme $H\overset{f_{|H}}{\longrightarrow}C$
induit par $f$ soit \`a fibres connexes (condition automatiquement v\'erifi\'ee).\\
Par la proposition 1.4, nous savons que $K_{X}^{2}\equiv0$ et donc,
puisque $D_{f}$ n'a que des composantes verticales, on a aussi
$D_{f}^{2}\equiv0$. On en d\'eduit que
$D_{|H}^{2}=\sum_{i}{D^{2}_{P_{i}}}_{|H}=0$
et, par le lemme de Zariski ([BPV]), que $D_{P_{i}}^{2}=0$ pour tout $i$. Par
ce m\^eme lemme de Zariski, on en d\'eduit alors qu'il existe
$r_{i}\in\mathbb{Q}^{*}$ tel que
${D_{P_{i}}}_{|H}=r_{i}({f_{|H}}^{*}P_{i})$. Par choix de la section
hyperplane, on obtient que $n_{i,j}=n_{i,k}=n_{i}$ pour tout triplet
$(i,j,k)$ et donc que $f^{*}P_{i}=n_{i}{(f^{*}P_{i})}_{red}$ et
$D_{P_{i}}=(n_{i}-1){(f^{*}P_{i})}_{red}$.\\
\indent Prenons un point $P_{i}\in C$ tel que la fibre sch\'ematique
$f^{*}P_{i}=n_{i}D_{i}$ soit non r\'eduite. Choisissons une boule
$\triangle$ contenant le point $P_{i}$ et posons
$U=f^{-1}(\triangle)$. On suppose que $P_{i}$ est le seul point de
$\triangle$ au dessus duquel la fibre est non r\'eduite. En appliquant
le lemme 4.3, on obtient alors un diagramme commutatif:  
\begin{equation*}
\begin{CD}
\overline{U} @)\overline{q}))U \\
@VV\overline{f}V @VVfV\\
\triangle @)z\mapsto z^{n_{i}}))\triangle
\end{CD}
\end{equation*}
\noindent o\`u $\overline{q}$ est un rev\^etement \'etale fini et
$\overline{f}$ est \`a fibres r\'eduites et connexes. Il n'est pas
difficile de voir que
$\overline{q}^{*}\omega_{U}=\overline{f}^{*}\omega_{\triangle}$. En
outre, puisque $\overline{q}$ est \'etale, l'injection
$\omega_{U}\hookrightarrow\Omega_{U}^{1}$ se rel\`eve en une 
injection 
$\overline{q}^{*}\omega_{U}=\omega_{\overline{U}}\hookrightarrow\overline{q}^{*}\Omega_{U}^{1}=\Omega_{\overline{U}}^{1}$
compatible \`a l'injection $\overline{f}^{*}\omega_{\triangle}\hookrightarrow\Omega_{\overline{U}}^{1}$.
On en d\'eduit alors que les fibres du morphisme $\overline{f}$ sont
normales (et donc irr\'eductibles), puisque 
les points singuliers des fibres de $\overline{f}$
correspondent aux points o\`u le conoyau de l'injection
$\overline{f}^{*}\omega_{\triangle}\hookrightarrow\Omega_{\overline{U}}^{1}$
n'est pas localement libre et puisque la structure de Poisson est
suppos\'ee quasi-r\'eguli\`ere. Il en r\'esulte que $F_{red}$ est une
surface normale, $F$ d\'esignant une fibre quelconque de $f$.\\
\indent En outre, si $F$
d\'esigne une fibre lisse de $f$, la formule 
d'adjonction fournit l'\'egalit\'e $\omega_{F}={\omega_{X}}_{|F}$ et donc
$\omega_{F}=\mathcal{O}_{F}$ puisque
$\mathcal{O}_{F}(D_{f})=\mathcal{O}_{F}$ dans ce cas. Par suite, $F$
est une surface $K3$ ou une surface ab\'elienne.\\
\indent Supposons maintenant qu'une fibre de $f$ soit une surface
$K3$. Par suite, pour toute fibre $F$ de $f$, on a
$\chi(\mathcal{O}_{F})=2$. Supposons que la fibre sch\'ematique $F$ soit
non r\'eduite. On peut donc \'ecrire $F=n(F_{red})$ avec $n\ge2$. A l'aide
du th\'eor\`eme de Riemann Roch, il n'est pas difficile de voir que
$\chi(\mathcal{O}_{F})=n\chi(\mathcal{O}_{F_{red}})$. En outre, par
la formule d'adjonction, on a
$\omega^{0}_{F_{red}}=\omega_{X}\otimes\mathcal{O}_{X}(F_{red})\otimes\mathcal{O}_{F_{red}}=f^{*}\omega_{C}(D_{f}+F_{red})\otimes\mathcal{O}_{F_{red}}=f^{*}\omega_{C}(n(F_{red}))\otimes\mathcal{O}_{F_{red}}=\mathcal{O}_{F_{red}}$
car $n(F_{red})=F$. Il en r\'esulte alors que
$\chi(\mathcal{O}_{F_{red}})=2$ par un th\'eor\`eme de Y.Umezu ([U] cor. 1) et
donc $n=1$. Aussi, le morphisme $f$ est non ramifi\'e et on a la formule
$\omega_{X}=f^{*}\omega_{C}$. Le morphisme $f$ est donc
lisse et analytiquement localement trivial ([F] thm. 4.8), ce que l'on voulait
d\'emontrer. En outre, la structure de Poisson est r\'eguli\`ere puisque,
comme nous l'avons d\'ej\`a remarqu\'e, les points o\`u la structure n'est pas
r\'eguli\`ere correspondent aux singularit\'es des fibres de $f$.\\
\indent Supposons maintenant qu'une fibre de $f$ soit une surface ab\'elienne.
Par un r\'esultat de J.Koll\'ar ([Ko]), le
faisceau $R^{1}f_{*}\mathcal{O}_{X}$ est localement libre de rang
2 et sa formation commute aux changements de base. Une
fibre r\'eduite $F$ de $f$ v\'erifie donc $h^{1}(F,\mathcal{O}_{F})=2$ et
est donc lisse par un r\'esultat de Y.Umezu ([U] cor. 1) puisque son faisceau
dualisant relatif est trivial. Si $F$ est une fibre de $f$ non
r\'eduite, on pose $F=n(F_{red})$ et $P=f(F)$. Choisissons une boule
$\triangle\subset C$ contenant le point $P$ et posons
$U=f^{-1}(\triangle)$. On suppose que $P$ est le seul point de
$\triangle$ au dessus duquel la fibre est non r\'eduite. En appliquant
le lemme 4.3, on obtient alors un diagramme commutatif:  
\begin{equation*}
\begin{CD}
\overline{U} @)\overline{q})) U \\
@VVf_{k}V @VVfV\\
\triangle @)z\mapsto z^{n_{i}})) \triangle
\end{CD}
\end{equation*}
\noindent tel que $\overline{q}$ soit un rev\^etement \'etale fini, et tel
que $\overline{f}$ soit \`a fibres r\'eduites et connexes. Il n'est pas
difficile de voir que
$\overline{q}^{*}\omega_{U}=\overline{f}^{*}\omega_{\triangle}$ puisque
nous avons la relation $\omega_{X}=f^{*}\omega_{C}(D_{f})$. Le
morphisme $\overline{q}$ \'etant \'etale on a
$\overline{q}^{*}\omega_{U}=\omega_{\overline{U}}$ et donc
$\omega_{\overline{U}}=\overline{f}^{*}\omega_{\triangle}$ de sorte
que les fibres de $\overline{f}$ ont toutes un faisceau dualisant
trivial. Enfin, puisque la formation du faisceau $R^{1}f_{*}\mathcal{O}_{X}$
commute aux changements de base, on en d\'eduit que toutes les fibres
$\overline{F}$ de $\overline{f}$ v\'erifient
$h^{1}(\overline{F},\mathcal{O}_{\overline{F}})=2$ ce qui implique
encore par le r\'esultat de Y.Umezu ([U] cor. 1) que ces fibres sont lisses. On
en d\'eduit alors que les fibres singuli\`eres de $f$ sont multiples de
surfaces lisses. Or, nous savons que les
singularit\'es des fibres du morphisme $f$ (munies de leurs structures
r\'eduites) correspondent aux singularit\'es de l'injection de faisceaux
$\omega_{X}\hookrightarrow\Omega^{1}_{X}$ (cf. corollaire 1.2), ce qui
prouve que la structure de Poisson est r\'eguli\`ere. Enfin, il n'est pas
difficile de voir que le fibr\'e canonique de $F_{red}$ est trivial pour
toute fibre $F$ de $f$ de sorte que, pour des raisons de caract\'eristique
d'Euler-Poincar\'e, les fibres de $f$ sont soit des surfaces ab\'eliennes,
soit des multiples de surfaces ab\'eliennes.\qed\\
\newline
\indent Les deux propositions qui suivent compl\`etent naturellement les
r\'esultats de la proposition pr\'ec\'edente.\\
\newline
\indent\textbf{Proosition 4.5}\textit{ Soit $X$ une vari\'et\'e
  projective de dimension 3. On suppose qu'il existe un
  morphisme surjectif
  $X\overset{f}{\longrightarrow}C$ vers une courbe alg\'ebrique $C$, dont les
  fibres sont soit des surfaces ab\'eliennes, soit multiples de surfaces
  ab\'eliennes. Enfin, on suppose qu'on a la formule
  $\omega_{X}=f^{*}\omega_{C}(D_{f})$, $D_{f}$ \'etant le diviseur de
  ramification de $f$.\\
\indent Alors $X$
est de la forme $(\overline{C}\times A)/G$ o\`u $\overline{C}$ est
une courbe lisse, $A$ une surface ab\'elienne
et $G\subset\text{Aut}(\overline{C})$ un groupe fini
op\'erant librement sur
$\overline{C}\times A$ par la formule:
$$g.(c,a)=(g.c,t_{g}(c,a)+u_{g}(a)),\,g\in G,\,c\in\overline{C},\,a\in A,$$
o\`u $u_{g}$ est un automorphisme de groupes de $A$ respectant la
structure symplectique et $t_{g}$ une fonction r\'eguli\`ere sur
$\overline{C}\times A$ \`a valeurs dans $A$.\\
\indent De plus, toute vari\'et\'e de cette forme admet une structure
de Poisson r\'eguli\`ere.}\\
\newline
\indent\textit{D\'emonstration}  Montrons que $X$ admet une
structure de Poisson 
r\'eguli\`ere. Consid\'erons l'injection de faisceaux
$\omega_{X}\cong f^{*}\omega_{C}(D_{f})\hookrightarrow\Omega_{X}^{1}$.
Le conoyau $\mathcal{G}$ de cette injection s'identifie, apr\`es le rev\^etement \'etale 
$\overline{X}\longrightarrow X$, au fibr\'e $\Omega^{1}_{\overline{X}}$
et il en r\'esulte que ce conoyau est localement libre de rang
2. Aussi, cette injection de fibr\'es vectoriels d\'efinira une structure
de Poisson r\'eguli\`ere si l'identit\'e de Jacobi est satisfaite. Par
le lemme 1.1, il suffit de v\'erifier l'int\'egrabilit\'e de
$\mathcal{G}^{*}\subset\mathcal{T}_{X}$ et donc, puisque $\overline{\pi}$
est \'etale, il suffit de v\'erifier que
$\overline{\pi}^{*}\mathcal{G}^{*}\subset\mathcal{T}_{\overline{X}}$
est int\'egrable, ce qui est imm\'ediat compte tenu de l'isomorphisme
$\overline{\pi}^{*}\mathcal{G}^{*}\cong\mathcal{T}_{\overline{X}/\overline{f}^{*}\omega_{\overline{C}}}$
compatible avec l'inclusion $\mathcal{T}_{\overline{X}/\overline{f}^{*}\omega_{\overline{C}}}\subset\mathcal{T}_{\overline{X}}$.\\
\indent Montrons maintenant qu'il existe un rev\^etement
galoisien ramifi\'e $\overline{C}\overset{\pi}{\longrightarrow}C$ tel que la
normalisation $\overline{X}$ du produit fibr\'e $\overline{C}\times_{C} X$
soit isomorphe au-dessus de 
$\overline{C}$ au produit $\overline{C}\times A$, $A$ \'etant une
surface ab\'elienne, et tel que la projection naturelle
$\overline{X}\longrightarrow X$ soit un rev\^etement \'etale.
Notons $P_{1},\ldots,P_{k}\in C$ les
points de $C$ tels que les fibres sch\'ematiques $f^{*}P_{i}$ soient non
r\'eduites. Soit $U=C\setminus\{P_{1},\ldots,P_{k}\}$ et soit
$X_{U}=f^{-1}(U)$. Par construction, le morphisme $f_{U}:=f_{|U}$ est
lisse. Notons $A\overset{a}{\longrightarrow}U$ le sch\'ema en groupes
au-dessus de la base $U$, dont les fibres sont les vari\'et\'es d'Albanese
des fibres correspondantes du morphisme $X_{U}\longrightarrow U$ et
notons $\mathcal{A}$ le faisceau de ses sections holomorphes. Alors,
le sch\'ema $X_{U}/U$ est un espace principal homog\`ene sous le sch\'ema en
groupes $A/U$ et correspond donc \`a un \'el\'ement $x$ du groupe de cohomologie
$H^{1}(U,\mathcal{A})$.\\
\indent  Montrons que cet \'el\'ement est de
torsion. Soit $W\subset U$ un ouvert de Zariski de $C$ inclus dans
$U$. Par le crit\`ere valuatif de propret\'e, la fl\`eche de restriction naturelle
$H^{1}(U,\mathcal{A})\overset{res}{\longrightarrow}H^{1}(W,\mathcal{A})$ est
injective. Le morphisme $X_{U}\overset{f_{U}}{\longrightarrow}U$ \'etant
lisse, il admet une quasi-section \'etale quasi-finie. Aussi, quitte \`a
restreindre l'ouvert $U$, on peut supposer qu'il existe une
section apr\`es rev\^etement \'etale connexe fini, de degr\'e
$m$. On en d\'eduit que l'\'el\'ement $x$ est de $m$-torsion.\\
\indent Notons $A_{m}/U$ le sous sch\'ema en groupes de $A/U$ form\'e des
\'el\'ements de $m$-torsion et $\mathcal{A}_{m}$ le faisceau des sections
correspondant. L'\'el\'ement $x\in H^{1}(U,\mathcal{A})$ provient donc
d'un \'el\'ement de $H^{1}(U,\mathcal{A}_{m})$.
Cet \'el\'ement d\'efinit un rev\^etement \'etale fini
$V\overset{\pi}{\longrightarrow}U$ qui trivialise l'espace principal
homog\`ene $X_{U}/U$. En outre, il existe une courbe lisse et connexe
$\overline{C}$ telle que le diagramme suivant soit commutatif:
\begin{equation*}
\begin{CD}
V @))) \overline{C}\\
@VV\pi V @VV\pi V\\
U @))) C
\end{CD}
\end{equation*}
o\`u les fl\`eches horizontales sont les inclusions naturelles.\\
\indent Il nous
reste \`a voir que le morphisme $\pi$ se ramifie convenablement aux
points $P_{i}$, \`a savoir que l'indice de ramification d'un point de
$\overline{C}$ au-dessus de $P_{i}$ est \'egal \`a la multiplicit\'e de la
fibre du morphisme $f$ en $P_{i}$, ce qui r\'esulte essentiellement du
lemme 4.3.\\
\indent Soit
 $\overline{X}$ la normalisation du produit 
fibr\'e $X\times_{\overline{C}}C$ et soit $\overline{f}$ le morphisme
naturel $\overline{X}\longrightarrow \overline{C}$. Alors,
$\overline{f}$ est un morphisme lisse et la projection naturelle 
$\overline{X}\overset{\overline{\pi}}{\longrightarrow}X$ est un
rev\^etement \'etale fini. Enfin, la formule
$\omega_{X}=f^{*}\omega_{C}(D_{f})$ entra\^{\i}ne l'\'egalit\'e
$\omega_{\overline{X}}=\overline{f}^{*}\omega_{\overline{C}}$ et le
morphisme $\overline{f}$ est donc localement trivial ([F] thm. 4.8). Mais, par
le crit\`ere valuatif de propret\'e, la fibration $\overline{X}/\overline{C}$ admet
une section et il en r\'esulte finalement que cette fibration est
triviale. La vari\'et\'e $\overline{X}$ est donc isomorphe, au-dessus de
$\overline{C}$ \`a $\overline{C}\times A$ o\`u $A$ est une surface
ab\'elienne. Remarquons enfin que le rev\^etement
$\overline{C}\longrightarrow C$ est galoisien et que son groupe $G$
agit naturellement sur $\overline{X}=\overline{C}\times A$ de mani\`ere
compatible \`a son action sur $\overline{C}$.\\
\indent Prenons $g\in G$ et $(c,a)\in\overline{C}\times
A$. Il n'est pas difficile de voir que
$g.(c,a)=(g.c,t_{g}(c,a)+u_{g}(a))$ o\`u $u_{g}$ est un automorphisme de
groupes et $t_{g}$ est une fonction
r\'eguli\`ere sur $\overline{C}\times A$ \`a valeurs dans $A$. Remarquons
alors que la 
structure de Poisson que nous avons construite sur $X$ se rel\`eve en
une structure de Poisson sur $\overline{X}$ qui est produit d'une
structure symplectique sur $A$ et de la structure nulle sur
$\overline{C}$ et que les automorphismes du rev\^etement
$\overline{X}\longrightarrow X$ sont des morphismes de Poisson. Un
calcul \'el\'ementaire montre alors que les \'el\'ements $u_{g}$ respectent la
structure symplectique sur $A$, ce qui est l'assertion souhait\'ee.\\
\indent  Enfin, l'assertion r\'eciproque est imm\'ediate.\qed\\
\newline
\indent\textbf{Proposition 4.6}\textit{ Soit $X$ une vari\'et\'e projective
  de dimension 3. On suppose qu'il existe un
  morphisme surjectif et lisse
  $X\overset{f}{\longrightarrow}C$ vers une courbe alg\'ebrique 
  $C$, dont les fibres sont des surfaces $K3$. On suppose en outre
  qu'on a la formule $\omega_{X}=f^{*}\omega_{C}$.\\
\indent Alors $X$ est de la forme $(\overline{C}\times S)/G$ o\`u
$\overline{C}$ est 
une courbe alg\'ebrique, $S$ est une surface $K3$ et $G$ un
groupe fini op\`erant librement sur $\overline{C}$, op\'erant sur $S$ en
respectant la structure symplectique et sur $\overline{C}\times S$ par
le produit de ses actions sur chacun des facteurs.\\
\indent De plus, toute vari\'et\'e de cette forme admet une structure
de Poisson r\'eguli\`ere.}\\
\newline
\indent\textit{D\'emonstration} Un argument analogue \`a celui utilis\'e
dans la proposition pr\'ec\'edente permet de montrer que $X$ admet une
structure de Poisson r\'eguli\`ere correspondant \`a l'injection de fibr\'es
vectoriels $\omega_{X}=f^{*}\omega_{C}\hookrightarrow \Omega_{X}^{1}$.\\ 
\indent Rappelons qu'il existe un espace de
modules fin $\mathcal{K}_{d,n}$ pour les surfaces $K3$ munies d'une
polarisation
de degr\'e $d$ fix\'e et d'une trivialisation du syst\`eme local
$H^{2}(.,\mathbb{Z}/n\mathbb{Z})$ avec $n\ge3$ fix\'e. Cela r\'esulte 
de l'existence d'un espace de modules fin pour les surfaces $K3$
marqu\'ees ([Be3] expos\'e VIII prop 1), du
fait qu'un automorphisme d'ordre fini induisant l'identit\'e sur
$H^{2}(.,\mathbb{Z}/n\mathbb{Z})$ est l'identit\'e par un lemme de Serre
([G], Appendice) et du fait que le groupe des automorphismes
projectifs d'une surface $K3$ est fini. Puisque le groupe
$\text{Aut}(H^{2}(.,\mathbb{Z}/n\mathbb{Z})$ est fini, il existe un
rev\^etement \'etale fini et connexe $\overline{C}\longrightarrow C$ de
$C$ tel que la famille $X\times_{C}\overline{C}/\overline{C}$ soit
rigidifi\'ee. Par suite, il existe un morphisme
$\overline{C}\longrightarrow\mathcal{K}_{d,n}$. Mais, puisque le
morphisme $f$ est
localement trivial ([F] thm. 4.8), ce morphisme est constant et, puisque
l'espace de modules est fin, la famille
$X\times_{C}\overline{C}/\overline{C}$ est triviale.\\
\indent On peut toujours supposer que le rev\^etement
$\overline{C}\longrightarrow C$ est galoisien de groupe $G$. On sait
que $\overline{C}\times_{C} X\cong \overline{C}\times S$ o\`u $S$ est une
surface $K3$. Le groupe $G$ agit sur le produit
 $\overline{C}\times S$ et sur la courbe $\overline{C}$ de mani\`ere
compatible. Puisque le groupe des automorphismes d'une surface $K3$
est discret, le groupe $G$ agit en fait sur $S$ et son action sur le
produit $S\times\overline{C}$ est finalement le produit des
actions sur chacun des facteurs. De plus, la structure de Poisson sur
$X$ induit une structure de Poisson r\'eguli\`ere sur
$S\times\overline{C}$. Il n'est pas difficile de voir que sur une
telle vari\'et\'e, une structure de Poisson r\'eguli\`ere est n\'ecessairement
le produit d'une structure symplectique sur $S$ et de la structure
triviale sur $\overline{C}$ et que le groupe $G$ respecte la
structure symplectique sur $S$, ce qui termine la preuve de notre
proposition, puisque la derni\`ere assertion est \'evidente.\qed\\ 
\newline
\indent Nous avons donc prouv\'e le th\'eor\`eme:\\
\newline
\indent\textbf{Th\'eor\`eme 4.7}\textit{ Soit $X$ une vari\'et\'e projective
 de dimension 3. On suppose en outre que $X$ est de dimension
  de Kodaira $\kappa(X)=1$.\\
\indent Alors $X$ admet une structure de Poisson quasi-r\'eguli\`ere si
et seulement si $X$ appartient \`a l'une des deux familles suivantes:\\
\indent (1) $X=(C\times S)/G$ o\`u $C$ est une courbe de genre au moins
2, $S$ est une surface $K3$ et $G$ un groupe fini op\'erant librement sur $C$
 et op\'erant sur $S$ en respectant la structure symplectique,\\
\indent (2) $X=(C\times A)/G$ o\`u $C$ est une courbe de genre au moins
2, $A$ est une surface ab\'elienne
et $G\subset\text{Aut}(C)$ un groupe fini op\`erant librement sur
$C\times A$ par la formule:
$$g.(c,a)=(g.c,t_{g}(c,a)+u_{g}(a)),\,g\in G,\,c\in C,\,a\in A,$$ 
o\`u $u_{g}$ est 
un automorphisme de groupes de $A$ respectant la structure
symplectique 
et  $t_{g}$ une fonction r\'eguli\`ere sur $C\times A$ \`a valeurs dans $A$.\\
\indent Enfin, $X$ est alors minimale, et la structure de Poisson est
r\'eguli\`ere.}\\
\newline
\indent\textbf{Remarque 4.8}  Les groupes finis d'automorphismes
symplectiques d'une surface $K3$ ont
\'et\'e classifi\'es par S.Mukai ([Mu]) et on dispose donc d'une
description assez pr\'ecise de ces vari\'et\'es.
\vspace{1cm}\\
\centerline{\large\textbf{R\'ef\'erences}}
$\ $
\newline
\noindent [BB] P.Baum, R.Bott, \emph{On the zeroes of meromorphic
  vector-fields}, Essays on Topology and related Topics, Springer,
29-45, 1970\\ 
\newline
[BPV] W.Barth, C.Peters, A.Van de Ven, \emph{Compact
  complex surfaces}, Springer-Verlag, 1984.\\
\newline
[Be1] A.Beauville, \emph{Surfaces alg\'ebriques complexes}, Ast\'erisque
54, 1978.\\
\newline
[Be2] A.Beauville, \emph{Vari\'et\'es de Prym et Jacobiennes
  interm\'ediaires}, Ann. scient. Ec. Norm. Sup., $4^{\textit{i\`eme}}$ s\'erie, 10,
309-391, 1977.\\
\newline
[Be3] A.Beauville, JP.Bourguignon, M.Demazure, \emph{G\'eom\'etrie des
  surfaces $K3$: modules et periodes}, Ast\'erisque 126, 1985.\\
\newline
[F] T.Fujita, \emph{On kähler fiber spaces over curves},
J. Math. Soc. Japan, 30, 779-794, 1978.\\
\newline
[G] A.Grothendieck, \emph{Construction de l'espace de Teichmuller},
S\'eminaire H.Cartan, 13, 1960-61.\\
\newline
[H] R.Hartshorne, \emph{Algebraic Geometry}, Graduate Text
  in Mathematics, Springer-Verlag, 1977.\\
\newline
[Ka1] Y.Kawamata, \emph{Minimal models and the Kodaira dimension of
  algebraic fiber spaces} J. de Crelle, 363, 1-46, 1985.\\
\newline
[Ka2] Y.Kawamata, \emph{Kodaira dimension of algebraic fiber spaces
  over curves}, Invent. Math., 66, 57-71, 1982.\\
\newline
[Ka3] Y.Kawamata, \emph{Pluricanonical systems on minimal algebraic
  varieties}, Invent. Math., 79, 567-588, 1985.\\
\newline
[Ko] J.Koll\'ar, \emph{Higher direct images of dualizing sheaves},
Ann. of Math., 123, 11-42, 1986.\\
\newline
[MP] Y.Miyaoka, T.Peternell, \emph{Geometry of Higher Dimensional
  Algebraic Varieties}, DMV seminar, 26, Birkauser, 1997.\\
\newline
[Mo] S.Mori, \emph{Threefolds whose canonical bundles are not
  numerically effective}, Ann. of Math., 116, 133-176, 1982.\\
\newline
[Mu] S.Mukai, \emph{Finite groups of automorphisms of K3 surfaces and
  the Mathieu group}, Invent. Math., 94, 183-221, 1988.\\
\newline
[Mum] D.Mumford, \emph{Abelian varieties}, Oxford Univ. Press, 1970.\\
\newline
[R] M.Reid, \emph{Bogomolov's theorem $c_{1}^{2}\le4c_{2}$},
Intl. Sympos. on Algebraic Geometry Kyoto, 623-642, 1977.\\
\newline
[U] Y.Umezu, \emph{On normal projective surfaces with trivial
  dualizing sheaf}, Tokyo J.Math, 4, 343-354, 1981.

\end{document}